\documentclass[graybox]{svmult}

\usepackage{type1cm}        
\usepackage{makeidx}         
\usepackage{graphicx}        
\usepackage{multicol}        
\usepackage[bottom]{footmisc}

\usepackage{newtxtext}       %
\usepackage[varvw]{newtxmath}       

\makeindex             

\newcommand{\RR}{\mathbb{R}}

\newcommand{\PP}{\mathbb{P}}
\newcommand{\BB}{\mathbb{B}}

\newcommand\FOCS{{\operatorname{\PP_1}}}
\newcommand\HOCS{{\operatorname{\BB_2}}}
\newcommand\AS{{\operatorname{AS}}}
\newcommand\SAS{{\operatorname{SAS}}}

\newcommand\SHS{{\operatorname{SHS}}}

\graphicspath{ {./Figures/} }
\usepackage{lipsum}
\usepackage[disable]{todonotes}
\usepackage{cleveref}

\begin{document}

\title*{Coarse Spaces Based on Higher-Order Interpolation for Schwarz Preconditioners for Helmholtz Problems}

\titlerunning{Higher-Order Coarse Spaces for Schwarz Preconditioners for Helmholtz Problems}

\author{Erik Sieburgh, Alexander Heinlein, Vandana Dwarka and Cornelis Vuik}
\institute{Erik Sieburgh, Alexander Heinlein, Vandana Dwarka and Cornelis Vuik \at Delft University of Technology, Delft Institute of Applied Mathematics, Mekelweg 4, 2628 CD Delft, Netherlands, \email{eriksieburgh@gmail.com, a.heinlein@tudelft.nl, v.n.s.r.dwarka@tudelft.nl, c.vuik@tudelft.nl}}
%
%
\maketitle

\abstract{The development of scalable and wavenumber-robust iterative solvers for Helmholtz problems is challenging but also relevant for various application fields. In this work, two-level Schwarz domain decomposition preconditioners are enhanced by coarse space constructed using higher-order B\'ezier interpolation. The numerical results indicate numerical scalability and robustness with respect the wavenumber, as long as the wavenumber times the element size of the coarse mesh is sufficiently low. 
}

\section{Introduction}\label{sec:introduction}
The Helmholtz equation, though seemingly simple, presents significant numerical challenges, particularly at large wavenumbers. These challenges can be attributed to two main problems. Firstly, a pollution error, a type of numerical dispersion error due to discrepancies between actual and numerical wavenumbers, necessitates grid refinement as the wavenumber increases, enlarging the linear system~\cite{Deraemaeker1999}. The second problem is the slow convergence behavior of iterative solvers for large wavenumbers and matrix dimensions. 

Despite extensive research into numerical solvers for Helmholtz problems, no scalable and wavenumber-robust solver has yet been developed; the situation becomes even more difficult when considering heterogeneities or the elastic Helmholtz equation. 

While the convergence of the conjugate gradient (CG) method for symmetric positive definite problems can be bounded via the condition number of the preconditioned matrix, $\kappa(M^{-1}A)$, the convergence is much less understood for solving indefinite problems using, for instance, the generalized minimal residual (GMRES) method~\cite{saad_gmres_1986}.

For GMRES, if the coefficient matrix $A$ is non-normal but diagonalizable, its convergence can resemble that of a normal matrix. Thus, the convergence rate cannot be directly bounded via the eigenvalues. However, they significantly influence it; cf.~\cite{DwaTwoLevel2020}.

Recently, in~\cite{DwaTwoLevel2020,DwaMulti2022}, higher-order Bézier deflation vectors haven been employed to design nearly wavenumber-independent multigrid solvers. This is facilitated by the use of a higher-order coarse interpolation. Moreover, there has also been work in Schwarz domain decomposition preconditioners for Helmholtz problems with high wavenumbers; see, for instance,~\cite{erlangga2008advances,graham2017domain,graham2017recent,Bonazzoli2018,Graham2020}. 

In this work, we employ coarse spaces for two-level Schwarz preconditioners spanned by higher-order Bézier functions. We observe good convergence properties for high wavenumber cases and increasing numbers of subdomains.

\section{Problem Description}\label{sec:problem_description}
\begin{figure}[t]
  \begin{center}
    \includegraphics[width=0.5\textwidth]{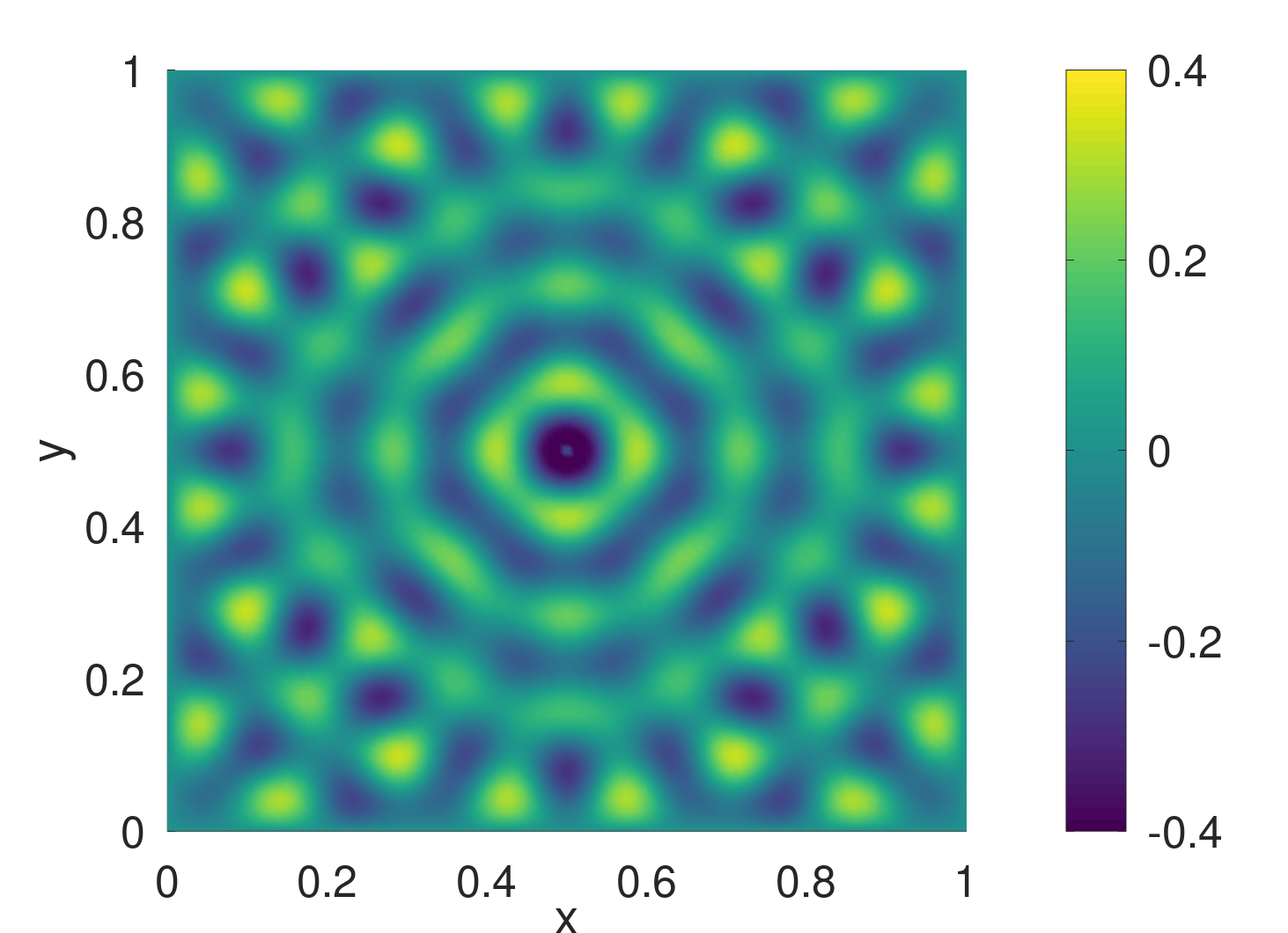}
  \end{center}
  \caption{Analytical solution of MP-1 with $k = 50$ on the unit square domain with homogeneous Dirichlet boundary conditions and the point source term located at the centre of the domain.}
  \label{fig:exact_sol_mp1}
\end{figure}

We consider the two-dimensional Helmholtz equations on the computational domain $\Omega = [0,1]\times[0,1] \subset \RR^2$, with constant wavenumber $k > 0 $ and either a Dirichlet boundary condition or a Sommerfeld radiation condition. As the first model problem, we seek the unknown field $u(x,y)$ that satisfies
\begin{eqnarray}\label{eq:BVP_MP-1}
    \left\lbrace
    \begin{array}{rll}
    -\nabla u(x,y)  -k^2 u(x,y) & = \delta(x - \frac{1}{2},y - \frac{1}{2}), & (x,y) \in \Omega, \\
    u(x,y) & = 0 &\text{ for } (x,y) \in \partial\Omega. 
    \end{array}
    \right.
\end{eqnarray}
A standard Dirac delta function is used as the point source term.
We refer to this boundary value problem (BVP) as \emph{MP-1}, and its analytical solution is depicted in~\cref{fig:exact_sol_mp1}.
For the second model problem, we replace the Dirichlet boundary condition with the Sommerfeld radiation condition given by
\begin{eqnarray*}\label{eq:BVP_MP-2}
\left(\frac{\partial u(x,y)}{\partial \textbf{n}} - i k u(x,y) \right) = 0,\quad \text{ for } (x,y)\in \partial\Omega,
\end{eqnarray*}
where $\textbf{n}$ is the outward pointing normal unit vector. We will refer to the model problem with Sommerfeld radiation condition as \emph{MP-2}.

We use a uniform structured grid $G_h$, where $h$ in $G_h$ denotes the mesh size $h = 1/(n-1)$; and $n$ is the number of grid nodes in one dimension. Then, we discretize using a second-order finite difference scheme, resulting in a linear system of equations
\begin{eqnarray}
    \label{eq:LS_final}
    A \Vec{u} = \Vec{f},
\end{eqnarray}

For both MP-1 and MP-2, the matrix $A$ is sparse and symmetric.
However, when using Dirichlet boundary conditions the linear system is indefinite real symmetric, while using the Sommerfeld radiation conditions results in an indefinite complex symmetric but non-Hermitian system. We consider both model problems because MP-2 resembles real-life infinite domain wave propagation more,
while MP-1 is more difficult for iterative solvers due to the absence of damping.

One major problem in numerically solving Helmholtz problems is the pollution error, which is a type of numerical dispersion error that accumulates due to discrepancies between the actual and numerical wavenumbers. 
For second-order finite difference discretizations, avoiding this error requires grid refinement with $k^3h^2 \leq 1$ \cite{Deraemaeker1999}. However, this leads to very large linear systems as $k$ increases. Therefore, the lighter condition $\kappa_h := kh \leq C$ is often used, implying a fixed number of grid nodes per wavelength $\lambda$.
The coarse grid in the two-level Schwarz preconditioners should also be sufficiently refined for increasing wavenumbers. We use $\kappa_H := kH = 1$~\cite{graham2017domain, graham2017recent}, with coarse grid size $H$.
Another way of reducing the pollution error is to adopt Isogeometric Analysis as a spatial discretization technique~\cite{Dwarka2021_IgA}.

\section{Two-Level Schwarz Preconditioners}\label{sec:domain_decomposition}


It has been shown that the classical one-level Schwarz domain decomposition methods are not scalable for Helmholtz problems; cf.~\cite{graham2017domain,graham2017recent}. 
In this research, we consider two-level Schwarz preconditioners which are constructed from a nonoverlapping domain decomposition 
$$
    \overline{\Omega} = \bigcup_{i=1}^N \overline{\Omega}_i
$$
as follows: On the first level, we extend the subdomains by layers of grid nodes to construct overlapping subdomains $\Omega_1,\ldots,\Omega_N$. Then, $R_i$ corresponds to the discrete restriction from $\Omega$ to $\Omega_i$, for $i=1,\ldots,N$. On the second level, we define a coarse space given by coarse basis functions being the rows of $R_0$, which are defined on a coarse mesh $G_{H}$ with the same step size as the subdomain size $H$. 

We consider different types of two-level overlapping Schwarz preconditioners, beginning with the two-level additive Schwarz (AS2) preconditioner given by
\begin{equation}\label{eq:DDM:SO_AS}
    M_{\mathrm{AS2}}^{-1} =  R_0^T (R_0A{R_0}^T)^{-1} R_0 + \displaystyle\sum_{i=1}^{N} R_i^T (R_iA{R_i}^T)^{-1} R_i.
\end{equation}
This corresponds to the case of exact local and coarse solvers, meaning that the local and the coarse problems are solved using a sparse direct solver.
For the Laplace problem and standard Lagrangian coarse basis functions, 
an upper bound for the condition number of the preconditioned system matrix is then given by
\begin{equation}\label{eq:CN_2LvL_AS}
    \kappa(M_{\mathrm{AS2}}^{-1}A) \leq C\left(1 + \frac{H}{\delta} \right);
\end{equation}
cf.~\cite[Theorem 3.13, p. 69]{Toselli2004}. 
This indicates numerical scalability, as the condition number is bounded by a constant when increasing the number of subdomains, as long as
the ratio of $H/\delta$ is kept fixed. Even though we do not have a theoretical proof, we will analyze the numerical scalability of two-level Schwarz preconditioners with the coarse spaces to be described in~\cref{sec:coarse_space} for Helmholtz problems in numerical experiments.

In addition to the AS2 preconditioner defined in~\cref{eq:DDM:SO_AS}, we also consider the two-level scaled additive Schwarz (SAS2) preconditioner
\begin{equation}\label{eq:DDM:SO_SAS}
    M_{\mathrm{SAS2}}^{-1} =  R_0^T (R_0A{R_0}^T)^{-1} R_0 + \displaystyle\sum_{i=1}^{N} R_i^T D_i (R_iA{R_i}^T)^{-1} R_i,
\end{equation}
with the diagonal matrices $D_i$ satisfying 
$\sum_{i=1}^{N} R_i^T D_i R_i = I$; cf.~\cite{cai_restricted_1999}. Furthermore, inspired by the use of deflation techniques for preconditioning Helmholtz problems in~\cite{DwaTwoLevel2020,Sheikh2013,graham2017domain,Graham2020},  we also consider the two-level scaled hybrid Schwarz (SHS2) preconditioner, defined as
\begin{align}\label{eq:DDM:SO_SHS}
    M_{\mathrm{SHS2}}^{-1} =\ & R_0^T (R_0A{R_0}^T)^{-1} R_0 + \nonumber \\
    &\left( \displaystyle \sum_{i=1}^{N} R_i^T D_i (R_iA{R_i}^T)^{-1}R_i \right) \left(I - A \left(R_0^T (R_0A{R_0}^T)^{-1} R_0\right) \right),
\end{align}
where the projection $P_0 := (I - A (R_0^T (R_0A{R_0}^T)^{-1} R_0) )$ deflates out the coarse space.
Expression
(\ref{eq:DDM:SO_SHS}) is a combination of~\cref{eq:DDM:SO_SAS} and the two-level hybrid Schwarz preconditioner from~\cite{Bootland2021}, which uses an ``adapted deflation technique''. The deflation should project the low frequency eigenvalues to zero and therefore remove the effect of those modes on the convergence of the iterative method.

\section{Coarse Space} \label{sec:coarse_space}

We define the coarse basis using higher-order B\'ezier interpolation on the coarse mesh taken from~\cite{DwaTwoLevel2020,DwaMulti2022} instead of classical linear coarse basis function; cf.~\cite{Quarteroni1999,Toselli2004}.
For the remainder of the article, the use of a linear coarse basis will be denoted as $\FOCS^H$, where $H$ is the coarse mesh size. 

\begin{figure}[h]
  \begin{center}
    \includegraphics[width=0.9\textwidth]{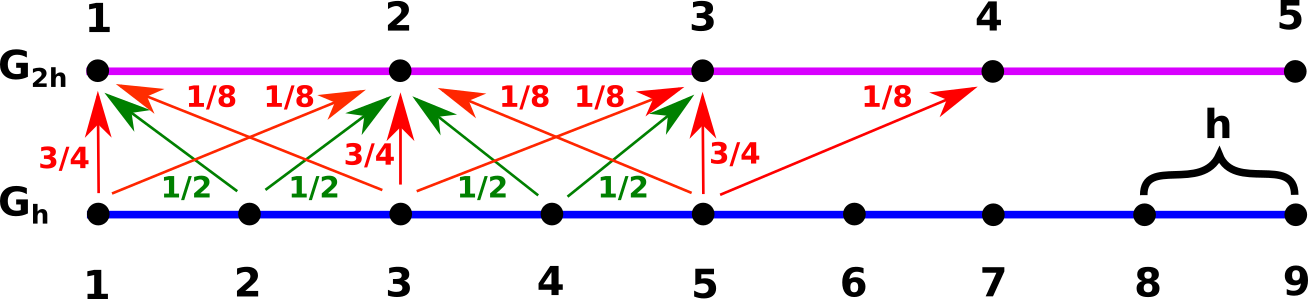}
  \end{center}
  \caption{Higher-order Bézier interpolation function for a fine and coarse grid with $H = 2h$.}
  \label{fig:HO_interpolation}
\end{figure}

The higher-order interpolation functions based on second-order rational B\'ezier curves from~\cite{DwaTwoLevel2020,DwaMulti2022}, for $H = 2h$, are given by
\begin{align}\label{eq:GCS:HO_R_0}
R_0\left[u_{h}\right]_{i}=\frac{1}{8}\Big(\left[u_{h}\right]_{(m-2)} &+  4\left[u_{h}\right]_{(m-1)} + \\ \nonumber
6\left[u_{h}\right]_{m} &+ 4\left[u_{h}\right]_{(m+1)} + \left[u_{h}\right]_{(m+2)} \Big),
\end{align}
for $i = 1,...,\frac{(n+1)}{2}$ and $m = 2i - 1$. This function for the restriction operator in~\cref{eq:GCS:HO_R_0} is partially visualized in~\cref{fig:HO_interpolation}. The prolongation from the coarse to the fine mesh is then given by the transpose $R_0^\top$.
The use of these higher-order Bézier interpolation functions in the preconditioners is denoted by $\HOCS^H$, where $H$ is again the coarse mesh size. Note that due to this implementation, the size of the coarse matrix is the same for both the linear and the higher-order Bézier interpolation function. However, the coarse matrix following from the higher-order Bézier interpolation function has a lower sparseness when compared to using the linear interpolation.

\section{Numerical Experiments}\label{sec:numerical_experiment}
In our numerical experiments, we test the performance of the two-level Schwarz preconditioners described in~\cref{sec:domain_decomposition} using the B\'ezier coarse interpolation defined in~\cref{sec:coarse_space} for the Helmholtz model problems MP-1 and MP2. We also compare the performance against using linear coarse basis functions.
Therefore, we employ the GMRES method with a relative tolerance of $10^{-7}$ and a maximum iteration cap of $100$ for the first two tables and $50$ for the others. These maximum iteration caps are on the low end, but since we are interested in wave-number independence having these maximum iteration caps is appropriate. 
In our tables, $\mathrm{x}$ indicates that the relative tolerance has not been reached within the maximum number of iterations. The implementation is done in MATLAB, and the computations have been carried out on a Intel Core i7-8557U processor with four CPU cores and 16\,GB memory.

\medskip

\noindent \textbf{Scalability of $\FOCS$ and $\HOCS$ coarse spaces}\quad For the numerical experiments, we use $H=4h$, $\kappa_h = 0.25$, and maximize the subdomain overlap $\delta$. The overlap is increased such that any grid node belongs to at most $4$ subdomains.

\begin{table}[t]
\centering
\begin{tabular}{r@{\hspace{5pt}}rcr@{\hspace{5pt}}r@{\hspace{5pt}}r@{\hspace{5pt}}r@{\hspace{5pt}}}
\noalign{\smallskip}\svhline\noalign{\smallskip}
\#\,subdomains & $|G_h|$ & \big/ & $|G_H|$ & $\frac{1}{h}$ & $\frac{1}{H}$ & $k$ \\
\noalign{\smallskip}\svhline\noalign{\smallskip}
400 & 6\,561 & \big/ & 441 & 80 & 20 &  20 \\
1\,600 & 25\,921 & \big/ & 1\,681 & 160 & 40 & 40 \\
3\,600 & 58\,081 & \big/ & 3\,721 & 240 & 60 & 60 \\
6\,400 & 103\,041 & \big/ & 6\,561 & 320 & 80 & 80 \\
10\,000 & 160\,801 & \big/ & 10\,201 & 400 & 100 & 100 \\
14\,400 & 231\,361 & \big/ & 14\,641 & 480 & 120 & 120 \\
19\,600 & 314\,721 & \big/ & 19\,881 & 560 & 140 & 140 \\
40\,000 & 641\,601 & \big/ & 40\,401 & 800 & 200 & 200 \\
\noalign{\smallskip}\hline\noalign{\smallskip}
\end{tabular}
\begin{tabular}{rrr}
\noalign{\smallskip}\svhline\noalign{\smallskip}
$\AS_{4h}$ & $\SAS_{4h}$ & $\SHS_{4h}$  \\
\noalign{\smallskip}\svhline\noalign{\smallskip}
21 & 15 & 9 \\
25 & 20 & 11\\
29 & 25 & 14\\
34 & 30 & 16\\
40 & 35 & 19 \\
45 & 41 & 22 \\
52 & 48 & 26 \\
75 & 71 & 38 \\
\noalign{\smallskip}\hline\noalign{\smallskip}
\end{tabular}
\begin{tabular}{rrr}
\noalign{\smallskip}\svhline\noalign{\smallskip}
 $\AS_{4h}$ & $\SAS_{4h}$ & $\SHS_{4h}$  \\
\noalign{\smallskip}\svhline\noalign{\smallskip}
 21 & 16 & 8 \\
 21 & 16 & 8\\
 21 & 16 & 8 \\
 21 & 16 & 8 \\
 21 & 16 & 8 \\
 21 & 16 & 8 \\
 21 & 16 & 8\\
 21 & 16 & 8 \\
\noalign{\smallskip}\hline\noalign{\smallskip}
\end{tabular}
\caption{Number of iterations results for preconditioners using $\FOCS$ (middle columns) or $\HOCS$ (right columns) interpolation for MP-2, with $\kappa_h = 0.25$, $H_{sub} = H = 4h$. $\mathrm{x}$ denotes the maximum number of iterations of $100$ being reached. Maximum subdomain overlap.}
\label{tab:MP2:kappa}
\end{table}
In~\cref{tab:MP2:kappa}, we compare iteration counts for MP-2 for linear coarse interpolation (middle columns) and higher-order B\'ezier grid interpolation (right columns).
As anticipated, for $\FOCS$, the number of iterations required to research convergence increases with the wavenumber, whereas
the convergence is robust with respect to the wavenumber when using higher-order B\'ezier grid interpolation. The iteration count is the lowest for the $\SHS_{4h}$ preconditioner, which requires to first solve the coarse problem before solving the local problems.

\begin{table}[t]
\centering
\begin{tabular}{r@{\hspace{5pt}}rcr@{\hspace{5pt}}r@{\hspace{5pt}}r@{\hspace{5pt}}r@{\hspace{5pt}}}
\noalign{\smallskip}\svhline\noalign{\smallskip}
\#\,subdomains & $|G_h|$ & \big/ & $|G_H|$ & $\frac{1}{h}$ & $\frac{1}{H}$ & $k$ \\
\noalign{\smallskip}\svhline\noalign{\smallskip}
400 & 6\,561 & \big/ & 441 & 80 & 20 &  20 \\
1\,600 & 25\,921 & \big/ & 1\,681 & 160 & 40 & 40 \\
3\,600 & 58\,081 & \big/ & 3\,721 & 240 & 60 & 60 \\
6\,400 & 103\,041 & \big/ & 6\,561 & 320 & 80 & 80 \\
10\,000 & 160\,801 & \big/ & 10\,201 & 400 & 100 & 100 \\
14\,400 & 231\,361 & \big/ & 14\,641 & 480 & 120 & 120 \\
19\,600 & 314\,721 & \big/ & 19\,881 & 560 & 140 & 140 \\
40\,000 & 641\,601 & \big/ & 40\,401 & 800 & 200 & 200 \\
\noalign{\smallskip}\hline\noalign{\smallskip}
\end{tabular}
\begin{tabular}{rrr}
\noalign{\smallskip}\svhline\noalign{\smallskip}
$\AS_{4h}$ & $\SAS_{4h}$ & $\SHS_{4h}$  \\
\noalign{\smallskip}\svhline\noalign{\smallskip}
 21 & 16 & 9 \\
 40 & 34 & 18\\
 55 & 49 & 26\\
 x & x & 53\\
 x & x & 55 \\
 x & x & 89 \\
 x & x & x \\
 x & x & x \\
\noalign{\smallskip}\hline\noalign{\smallskip}
\end{tabular}
\begin{tabular}{rrr}
\noalign{\smallskip}\svhline\noalign{\smallskip}
$\AS_{4h}$ & $\SAS_{4h}$ & $\SHS_{4h}$  \\
\noalign{\smallskip}\svhline\noalign{\smallskip}
 20 & 14 & 7 \\
 18 & 13 & 6\\
 19 & 13 & 7 \\
15 & 12 & 5 \\
 18 & 13 & 6 \\
 19 & 13 & 7 \\
 18 & 12 & 6\\
 19 & 13 & 7 \\
\noalign{\smallskip}\hline\noalign{\smallskip}
\end{tabular}
\caption{Number of iterations results for preconditioners using $\FOCS$ (middle columns) or $\HOCS$ (right columns) interpolation for MP-1, with $\kappa_h = 0.25$, $H_{sub} = H = 4h$. $\mathrm{x}$ denotes the maximum number of iterations of $100$ being reached. Maximum subdomain overlap.}
\label{tab:MP1:kappa}
\end{table}

Next, in~\cref{tab:MP1:kappa}, the corresponding results for 
MP-1 are listed. The results are similar to those for MP-2 in that the iteration count for the $\FOCS$ coarse space increases with the wavenumber, and the convergence for $\HOCS$ coarse space is robust. 


The subdomain sizes are relatively low, and hence, the dimension of the coarse space, $|G_H|$, 
is relatively high compared with the original problem size, $|G_h|$. Hence, in parallel computations, the solution of the coarse problem would quickly become a bottleneck. In distributed-memory parallel computations, we would typically have much more memory available on each rank, such that we could significantly increase the subdomain size, resulting in a better ratio of the coarse and fine grid sizes.


\medskip

\noindent \textbf{Detailed variation of $k$ and $n$}\quad Finally, we further analyze the performance of the best-performing preconditioner, $\SHS_{H}/\HOCS^{H}$ by varying the wavenumber $k$ and the grid size $n$. In particular, we investigate the limit of the wavenumber robustness of the preconditioner for MP-1. We vary the coarsening ratio $H/h$ between 4 and 16, as shown in \cref{tab:MP1_kappa_analysis_4h,tab:MP1_kappa_analysis_16h}. As expected, if $n$ is too low, the iteration count increases with the wavenumber. This is due to $h$ and $H$ being too large to maintain suitable $\kappa_h$ and $\kappa_H$ ratios.

Looking at ~\cref{tab:MP1_kappa_analysis_4h}, we observe that for $H = 4h$, $\kappa_h = 0.25$ seems to be suitable. From~\cref{tab:MP1_kappa_analysis_16h} on the other hand,
with $H = 16h$, $\kappa_h = 0.0625$ seems to be required for wavenumber robustness, which is $1/4$ of that for $H = 4h$. This seems to indicate that the wavenumber robustness also requires sufficient resolution of the coarse mesh rather than only for the fine mesh. In particular, the same $\kappa_H := kH \leq 1$ seems to be required for both subdomain sizes.

\begin{table}[t]
\centering
    \small \setlength{\tabcolsep}{3pt}
    \begin{tabular}{r|rrrrrrrrrrrrrrrrr}
    \hline $k \setminus n$ & $33 $ & $41$ & $49$ & $57$ & $65$ & $73$ & $81$ & $89$ & $97$ & $105$ & $113$ & $121$ & $129$ & $137$ & $145$ & $153$ & $161$\\
    \hline 
        10 & 7  & 6  & 6  & 6  & 6  & 6  & 6  & 6  & 5  & 5  & 5  & 5  & 5  & 5  & 5  & 5  & 5  \\
        20 & 18  & 17  & 11  & 9  & 8  & 7  & 7  & 7  & 7  & 7  & 7  & 6  & 6  & 6  & 6  & 6  & 6  \\
        30 & 38  & 23  & 32  & 37  & 17  & 11  & 9  & 8  & 7  & 7  & 7  & 7  & 6  & 6  & 6  & 6  & 6  \\
        40 & x & x & 45 & x & x & x & 43 & 15  & 10  & 9  & 8  & 7  & 7  & 7  & 7  & 6  & 6  \\
        50 & x & x & x & x & 48 & x & x & x & x & 31 & 14 & 11 & 9 & 9 & 8 & 8 & 7 \\
        \hline
    \end{tabular}
    \caption{Number of iterations with the $\SHS_{4h}/\HOCS^{4h}$ preconditioner for MP-1 for small $k$. $\mathrm{x}$ denotes the maximum NOI of $50$ being reached. Maximum subdomain overlap.}
    \label{tab:MP1_kappa_analysis_4h}
\end{table}

\begin{table}[t]
\centering
    \small \setlength{\tabcolsep}{3pt}
    \begin{tabular}{r|rrrrrrrrrrrrrrr}
    \hline $k \setminus n$ & $33$ & $49$ & $65$ & $81$ & $97$ & $113$ & $129$ & $145$ & $161$ & $177$ & $193$ & $209$ & $225$ & $241$ & $257$ \\
    \hline 
        5  & 5  & 8  & 7  & 8  & 7  & 8  & 7  & 8  & 7  & 8  & 7  & 8  & 7  & 8  & 7  \\
        10 & 6  & 11 & 10 & 13 & 9  & 8  & 7  & 7  & 7  & 7  & 7  & 7  & 7  & 7  & 7  \\
        15 & 11 & 11 & 12 & 17 & 12 & 16 & 12 & 12 & 9  & 10 & 9  & 9  & 8  & 9  & 8  \\
        20 & 7  & 13 & 14 & 18 & 16 & 20 & 18 & x & 21 & 17 & 12 & 11 & 9  & 10 & 9  \\
        25 & 8  & 16 & 21 & 26 & 18 & 45 & 19 & 28 & 24 & 31 & 24 & 23 & 14 & 13 & 10 \\
        30 & 9  & 22 & 33 & 28 & 39 & 38 & 49 & 36 & 33 & x & 34 & 40 & 40 & 29 & 18 \\
        \hline
    \end{tabular}
    \caption{Number of iterations with the $\SHS_{16h}/\HOCS^{16h}$ preconditioner for MP-1 for small $k$. $\mathrm{x}$ denotes the maximum NOI of $50$ being reached. Maximum subdomain overlap.}
    \label{tab:MP1_kappa_analysis_16h}
\end{table}

\section{Conclusion}\label{sec:conclusion}
We have shown that using higher-order B\'ezier interpolation for the coarse space in two-level Schwarz preconditioners can yield scalable and wavenumber robust coarse spaces for Helmholtz problems. In particular, wavenumber robustness can be obtained for cases where a linear coarse space does not yield robust convergence; nonetheless, we observed that a certain refinement level of the coarse mesh relative to the wavenumber is required. 

The best-performing preconditioner is the two-level scaled hybrid Schwarz preconditioner, which uses a higher-order approximation scheme to construct the deflation vectors for the coarse space.

%
%

\bibliographystyle{spmpsci}
\bibliography{references}

\end{document}